\documentclass[10pt]{article}

\usepackage[utf8]{inputenc}
\usepackage[T1]{fontenc}
\usepackage{amsmath, amsfonts, amssymb, amsthm}
\usepackage{stmaryrd}
\usepackage{bbold}
\usepackage{graphicx}
\usepackage[export]{adjustbox}
\usepackage{hyperref}
\usepackage{cite}
\usepackage{geometry}
\usepackage{babel}

\title{Critical Exponent Elliptic Equations on the Half-Space: Uniqueness and Explicit Solutions}

\author{Azam\ Nouri \\[4pt]
\small Department of Science, Technology \& Mathematics, Lincoln University}
\date{}

\begin{document}
\maketitle

\begin{abstract} 
\noindent 
We prove that all positive solutions of $-\Delta u=u^{(2 n)/(n-2)}$
on the upper half space $\mathbb{R}_{+}^{n}(n \geq 3)$ satisfying the boundary condition $D_{x_{n}}(u)=-u^{n /(n-2)}$ are of the form 
\[
u = a \left(\frac{\lambda}{{\lambda^{2}+|x-y|^{2}}}\right)^{(n-2)/2},
\]
where $a = a(n)$, $\lambda>0$ and $y = (y_1,\ldots,y_n)$ is some point in the lower half-space with $y_{n}<0$.
\end{abstract}

\subsection*{1. Introduction} 
Consider the problem
\begin{align}
\begin{cases}
-\Delta u=u^{\frac{2 n}{n-2}} & \text { in } \mathbb{R}_{+}^{n}(n \geq 3), \\ 
D_{x_{n}}(u)=-u^{p} & \text { on } \partial \mathbb{R}_{+}^{n}(p = n /(n-2)), \\ 
u>0
\end{cases}\tag{1.1}
\end{align}
where $\mathbb{R}_{+}^{n}=\left\{(x_{1}, \ldots, x_{n-1}, x_{n}) : x_{n} \geq 0\right\}$ denotes the upper half space and $\partial \mathbb{R}_{+}^{n}$ its boundary.\\

Equation (1.1) is a special case of the Lane-Emden equation which arises in the study of stellar structures, plasma physics and fluid dynamics. The equation also appears in the context of the Sobolev embedding theorem including the trace embedding theorem and in the context of Riemannian geometry.

Our goal is to show that any positive solution of equation (1.1) can be expressed in the form of (1.2), where $a = a(n)$, $\lambda>0$ and $y = (y_1,\ldots,y_n)$ is some point in the lower half-space with $y_{n}<0$.
\begin{align}
u = a \left(\frac{\lambda}{{\lambda^{2}+|x-y|^{2}}}\right)^{(n-2)/2}\tag{1.2}
\end{align}
We restrict our attention to $C^{2}$-solutions.

\paragraph{Relation to Ou \cite{ou1996}.}
This paper is a continuation of the work of Ou \cite{ou1996}. In that paper, Ou studied the \emph{harmonic} case
\[
-\Delta u = 0 \quad \text{in } \mathbb{R}_+^n,
\]
under the same nonlinear boundary condition $D_{x_n}u = -u^{p}$ on $\partial \mathbb{R}_+^n$. He proved that when $p=n/(n-2)$, all positive solutions are precisely fundamental solutions of the Laplace equation multiplied by constants, and that no positive solutions exist in the subcritical case $p<n/(n-2)$.  

Our present work extends this framework to the \emph{critical semilinear equation}
\[
-\Delta u = u^{\frac{2n}{n-2}} \quad \text{in } \mathbb{R}_+^n,
\]
with the same nonlinear boundary condition. The structure of the arguments—Kelvin transform, decay at infinity, moving planes, and axial symmetry—parallels that of Ou, and many of the proofs of lemmas and theorems are straightforward adaptations of his. To avoid redundancy, we omit those details and refer the reader to \cite{ou1996}.

\paragraph{Methodology.}
The overall methodology follows \cite{NouriDissertation2023, ou1996}, with the moving-plane method as the central tool. This technique, introduced by Gidas, Ni, and Nirenberg \cite{gidas1979} and refined in many subsequent works \cite{caffarelli1989,chen1991,li1993,li1995}, establishes symmetry and monotonicity of solutions to semilinear elliptic equations.  

Our analysis begins by applying Kelvin's transform to a solution of (1.1) with respect to a boundary point on $\partial\mathbb{R}_+^n$. The transformed function is positive, decays at the rate $|x|^{2-n}$ at infinity, and satisfies a nonlinear boundary condition except at a possible singularity.  

We then employ the moving-plane method on the transformed problem. For the critical exponent $p=n/(n-2)$, this shows that the transformed solution must be symmetric about an axis parallel to the $x_n$-axis. Reversing the Kelvin transform yields the same symmetry for the original solution of (1.1).  

Finally, by combining this symmetry with a conformal mapping argument and the maximum principle, we deduce that every solution must take the explicit form (1.2).

\subsection*{2. Kelvin's transform}
 The Kelvin's transform is defined as
 \begin{align}
v(x) = \frac{1}{|x|^{n-2}}u\left(\frac{x}{|x|^2}\right).
\tag{2.1}
\end{align}
 where $u$ is a solution of (1.1) and $v$ is the Kelvin transform of $u$. One can check $v$ satisfies \cite{ou1996}:
 \begin{align}
\begin{cases}
-\Delta v = v^{\frac{2n}{n-2}} & \text{in } \mathbb{R}_{+}^{n}, \\
D_{x_{n}}(v) = -\frac{1}{|x|^{n-p(n-2)}} v^{p} & \text{on } \partial \mathbb{R}_{+}^{n}, \\
v > 0
\end{cases}
\tag{2.2}
\end{align}
except that $v$ may be singular at the origin.

We consider a family of parallel planes perpendicular to $\partial\mathbb{R}_+^n$, and choose the coordinate system so that the equations of the planes are of the form $x_1=\lambda$. We denote $x^{\lambda}=\left(2\lambda - x_1, x_2,\ldots,x_n\right)$ as the reflection of $x$ about the plane $x_1=\lambda$. We also denote $v_{\lambda}(x)=v\left(x^{\lambda}\right)$ as the reflection of $v$ about the plane $x_1=\lambda$.
We define the domain $\Sigma_{\lambda}$ as
$$
\Sigma_{\lambda}=\left\{x \in \mathbb{R}_+^n : x_1 \geq \lambda, x\neq 0\right\}.
$$
By reflecting $v$ about the planes $x_1=\lambda$ for $\lambda < 0$, we obtain a family of functions $\{v_\lambda\}_{\lambda< 0}$, which satisfies the same PDE as $v$ but with a different boundary condition. We can use the moving-plane method \cite{gidas1979,caffarelli1989} to study the properties of the family $\{v_\lambda\}_{\lambda < 0}$ to show that $v_\lambda$ is symmetric about the plane $x_1=\lambda$.\\

\noindent\textbf{Theorem 2.3.} (Ou \cite{ou1996})  
(a) There exists a sufficiently large constant $N$ such that the inequality $v(x) \geq v_{\lambda}(x)$ holds on $\Sigma_{\lambda}$ for all $\lambda \leq -N$.  

(b) If there exists some $\lambda_{0} < 0$ such that $v(x) \geq v_{\lambda_{0}}(x)$ on $\Sigma_{\lambda_{0}}$, but $v(x) \not\equiv v_{\lambda_{0}}(x)$ on $\Sigma_{\lambda_{0}}$, then there exists an $\epsilon > 0$ such that $v(x) \geq v_{\lambda}(x)$ on $\Sigma_{\lambda}$ for all $\lambda$ in the interval $[\lambda_{0}, \lambda_{0} + \epsilon)$.  

In other words, (a) states that $v$ is greater than or equal to its reflections $v_{\lambda}$ for $\left|\lambda\right|$ sufficiently large, while (b) states that if $v$ is greater than or equal to one of its reflections $v_{\lambda_{0}}$, but not identically equal to it, then $v$ is greater than or equal to all the reflections $v_{\lambda}$ in a small interval of $\lambda$ around $\lambda_{0}$.  

\noindent\textbf{Proof.} See Ou \cite{ou1996}.\\

What we learn from Theorem 2.3 brings us to the following conclusions: 

\noindent\textbf{Case 1:} Assume that $u(x)$ is symmetric about any axis parallel to the $x_n$-axis. This means that $u(x)$ does not depend on the variables $x_1, x_2, \ldots, x_{n-1}$, but only on the variable $x_n$. It follows that $u(x)$ can be expressed as a one-dimensional function of the form $u(x) = u(x_n)$, which reduces to the following ODE:
$$
-u''(x_n) = \left(u(x_n)\right)^{\frac{2n}{n-2}}, \qquad 
u'(0) = -u^{p}(0).
$$
This ODE does not have a solution.\\

\noindent\textbf{Case 2:} Assume that $u(x)$ is symmetric about an axis parallel to the $x_n$-axis. We consider the function 
\[
v(x) = \frac{1}{|x|^{n-2}} u\left(\frac{x}{|x|^2}\right),
\]
obtained through Kelvin's transform. By applying the moving-plane method \cite{gidas1979,li1993}, we have established the existence of a value $\lambda_0 < 0$ such that $v_{\lambda_0}(x) = v(x)$.

Suppose, for the sake of contradiction, that $v(x)$ is symmetric about the $x_n$-axis. This would imply that $u(x)$ is also symmetric about the $x_n$-axis, contradicting our initial assumption. Hence, $v(x) = v_{\lambda_0}(x)$ for some $\lambda_0 < 0$.

Now, if $v_{\lambda_0}(x) = v(x)$ holds for some $\lambda_0 < 0$, it means that $v(x)$ remains the same under the reflection across the plane $x_1 = \lambda_0$. This indicates that the origin is a removable singularity.

In this case, $v$ does not have any singularity, and $v$ is axially symmetric about an axis parallel to the $x_{n}$-axis. Moreover, the Kelvin transform of $u$ with respect to any point on $\partial \mathbb{R}_{+}^{n}$ is also symmetric about an axis parallel to the $x_n$-axis.\\

In conclusion, when $p=n /(n-2)$, the function $v$ is symmetric about an axis parallel to the $x_{n}$-axis. The same proof implies that the Kelvin transform of $u$ with respect to any point on $\partial \mathbb{R}_{+}^{n}$ is symmetric about an axis parallel to the $x_{n}$-axis.

\subsection*{3. Solutions of $\mathbf{(1.1)}$}
In this section, we derive the specific form of $u$ on $\partial \mathbb{R}_{+}^{n}$.

\noindent\textbf{Lemma 3.1.}  
There is a positive constant $\mu>0$ such that the following limits hold:
\begin{align}
\begin{aligned}
& \lim _{x \rightarrow \infty}|x|^{n-2} u(x)=\mu, \\
& \lim _{x \rightarrow \infty}|x|^{n-2} D_{x_{i}} u(x)=0, \quad i=1, \ldots, n, \\
& \lim _{x \rightarrow \infty}|x|^{n-2} x \cdot \nabla u(x)=-(n-2) \mu .
\end{aligned}
\tag{3.1}
\end{align}
Furthermore, $u$ is symmetric about an axis parallel to the $x_{n}$-axis.\\

\noindent\textbf{Proof.} See Ou \cite{ou1996}.\\

From this point forward, let us assume that $u$ is symmetric about the $x_n$-axis. By direct computation, we observe the following relationship:
$$
u_{s}(x)=s^{(n-2) / 2} u(s x), \quad s>0.
$$
This equation expresses the transformation of the function $u$ under scaling by a factor of $s$.  
All functions $u_s$ solve equation (1.1) for $p=n/(n-2)$ and are symmetric about the $x_n$-axis.

Furthermore, if we consider the Kelvin transform of $u$ with respect to the point $e=(1,0,\ldots,0)$, denoted $v_s(x)$, we have
\begin{align}
v_{s}(x) = \frac{1}{|x-e|^{n-2}} u_{s}\left(e+\frac{x-e}{|x-e|^{2}}\right).
\tag{3.2}
\end{align}
The functions $v_s(x)$, defined as the Kelvin transform of $u_s$ with respect to $e$, are also solutions of (1.1) with $p=n/(n-2)$. These solutions are symmetric about axes that pass through points on the $x_1$-axis.\\

\noindent\textbf{Lemma 3.2.}  
There exists an $s>0$ such that $v_s$ is symmetric about the axis parallel to the $x_n$-axis and passing through the point $(0.5,0,\ldots,0)$.\\

\noindent\textbf{Proof.}
Consider the point $x_t=(t,0,\ldots,0)$. We define the one-dimensional functions $f(t)$ and $g_s(t)$ as functions generated by $u$ and $v_s$, respectively. 
$$
f(t)=u\left(x_{t}\right), \quad g_{s}(t)=\frac{s^{(n-2) / 2}}{|t-1|^{n-2}} f\left(s\left(1+\frac{1}{t-1}\right)\right)
$$
The monotonicity property associated with the moving-plane method allows us to focus on showing the existence of an $s>0$ such that $g_s'(0.5)=0$.

By elementary calculations, we obtain:
$$
\begin{aligned}
g_{s}^{\prime}(0.5) & =2^{n} s^{(n-2) / 2}\left(\frac{n-2}{2} f(-s)-s f^{\prime}(-s)\right)=2^{n} s^{(n-2) / 2}\left(\frac{n-2}{2} f(s)+s f^{\prime}(s)\right) \\
& =2^{n} s^{(n-2) / 2}\left(\frac{n-2}{2} u\left(x_{s}\right)+x_{s} \cdot \nabla u\left(x_{s}\right)\right) .
\end{aligned}
$$
We made use of the fact that $f(-s) = f(s) = u(x_s)$. It is evident that $g_s'(0.5)$ is positive if $s$ is sufficiently small. From equation (3.1), we can deduce that $g_s'(0.5)$ becomes negative if $s$ is sufficiently large. Therefore, by the Intermediate Value Theorem, there exists an $s$ for which $g_s'(0.5) = 0$.\\

\noindent\textbf{Theorem 3.3.}  
For some $a=a(n), \lambda>0, x_{0}=\left(x_{0}^{\prime},0\right) \in \mathbb{R}^{n}$, we have
\begin{align}
u(x', 0) = a\left(\frac{\lambda}{\lambda^2 + |x' - x_0'|^2}\right)^{(n-2)/2}, \quad \forall x' \in \mathbb{R}^{n-1}.
\tag{3.3}
\end{align}

\noindent\textbf{Proof.}  
By Lemma 3.2, we choose $s$ such that $v_s(x)$ is symmetric about the axis passing through $(0.5,0,\ldots,0)$ and parallel to the $x_n$-axis.

Recall that an inversion is a conformal mapping that preserves angles and maps circles or spheres to other circles or spheres. In particular, consider the point $x = \left(0.5, x_2, \ldots, x_{n-1}, 0\right)$. It can be verified directly that
$$
\left|e + \frac{x - e}{|x - e|^2}\right| = 1
$$
where $e=(1,0,\ldots,0)$.  

Thus, because of (3.2),
$$
v_{s}\left(0.5, x_{2}, \ldots, x_{n-1}, 0\right)=\frac{1}{\left(0.25+x_{2}^{2}+\cdots+x_{n-1}^{2}\right)^{(n-2) / 2}} s^{(n-2) / 2} f(s).
$$ 
By the symmetry of $v_{s}$,
$$
v_{s}\left(x_{1}, x_{2}, \cdots, x_{n-1}, 0\right)
=\frac{1}{\left(0.25+\left(x_{1}-0.5\right)^{2}+x_{2}^{2}+\cdots+x_{n-1}^{2}\right)^{(n-2) / 2}} s^{(n-2) / 2} f(s).
$$

Since $v_{s}$ is a solution of (1.1), by applying the inverse of Kelvin's transform, we obtain the desired result.\\

\paragraph{Analyticity and Uniqueness.}
Given an elliptic partial differential equation with (real) analytic coefficients, it is known that solutions to such equations are (real) analytic. A real analytic function is one that can be expressed as a power series expansion around any point in its domain, which is guaranteed to converge to the function in some neighborhood of the point.

The "analyticity of the solutions of elliptic equations" has been established through the seminal contributions of \cite{morrey1957,morrey1958,morrey1958b,petrowsky1939}.
Therefore, the solution $u$ that we obtained on the boundary $\partial \mathbb{R}^{n}_+$ in Theorem~3.3 is a real analytic function.

The Cauchy–Kowalevski theorem states that if the coefficients of a system of partial differential equations are analytic functions, and if the initial data are also given as analytic functions, then there exists a unique analytic solution of the system in some neighborhood of the initial surface. The theorem and its proof are valid for analytic functions of either real or complex variables. The Cauchy–Kowalevski Theorem is valid also for fully nonlinear analytic PDE (see \cite{folland2021}).

Since both the equation and the boundary condition in (1.1) have real analytic coefficients, we can apply the Cauchy–Kowalevski theorem to conclude the uniqueness of solution (3.3) in a neighborhood of any point on the boundary $\partial \mathbb{R}_{+}^{n}$. The theorem ensures that the solution can be extended to a larger domain where it remains real analytic. 

Note that if $f$ and $g$ are two real analytic functions that agree on an open subset $\Omega$ of $\mathbb{R}^{n}$, then they must agree on the entire space $\mathbb{R}^{n}$. For a detailed explanation on this topic, we refer the reader to \cite{krantz2002}.  

Thus, the global solution for the system (1.1) must be of form (1.2).  

\begin{center}
\textbf{Acknowledgment}
\end{center}
I want to express my gratitude to my advisor Dr. Biao Ou, for his guidance and support throughout the course of this research. His feedback and mentorship were crucial in the completion of this work.

\bibliographystyle{plain}

\end{document}